\documentclass{siamart171218}
\usepackage{graphicx}
\usepackage{epstopdf}
\usepackage{amssymb}
\usepackage{amsmath}
\usepackage{hyperref}
\usepackage{cite}

\numberwithin{theorem}{section}
\newtheorem{example}{Example}
\newtheorem{remark}{Remark}

\newcommand{\TheTitle}{A discrete Gr\"{o}nwall inequality with application to
numerical schemes for subdiffusion problems}

\title{{\TheTitle}\thanks{Submitted to the editors DATE.
\funding{This work was funded by a grant 1008-56SYAH18037
from NUAA Scientific Research Starting Fund of Introduced Talent and a grant DRA2015518
from 333 High-level Personal Training Project of Jiangsu Province;
Australian Research Council grant DP140101193;
NSFC grants 11771035, 91430216, U1530401.}}}
\author{
Hong-lin Liao\thanks{Department of Mathematics,
Nanjing University of Aeronautics and Astronautics,
Nanjing, 211106, P. R. China.
(\email{liaohl@csrc.ac.cn}).}
\and
William McLean\thanks{School of Mathematics and Statistics,
University of New South Wales, Sydney 2052, Australia.
(\email{w.mclean@unsw.edu.au}).}
\and
Jiwei Zhang\thanks{Beijing Computational Science Research Center,
Beijing, 100094, P. R. China.
(\email{jwzhang@csrc.ac.cn}).}}

\newcommand{\fd}[1]{\mathcal{D}^{#1}_t}
\newcommand{\dfd}[1]{\mathcal{D}^{#1}_\tau}
\newcommand{\ffd}[1]{D^{#1}_f}
\newcommand{\diff}{\triangledown_\tau}
\newcommand{\Ass}[1]{\textbf{\upshape A#1}}
\newcommand{\R}{\mathbb{R}}
\newcommand{\iprod}[1]{\left\langle#1\right\rangle}

\newcommand{\rhomax}{\rho}
\newcommand{\opL}{\mathcal{L}}
\newcommand{\Bilin}{\mathcal{B}}
\newcommand{\defeq}{:=}

\newcommand{\zd}{\,\mathrm{d}}
\newcommand{\iprodb}[1]{\big\langle#1\big\rangle}

\newcommand{\mynorm}[1]{\left\|#1\right\|}
\newcommand{\mynormb}[1]{\big\|#1\big\|}
\newcommand{\mynormt}[1]{\|#1\|}
\newcommand{\bra}[1]{\left(#1\right)}
\newcommand{\brab}[1]{\big(#1\big)}
\newcommand{\braB}[1]{\Big(#1\Big)}
\newcommand{\brat}[1]{(#1)}
\newcommand{\kbra}[1]{\left[#1\right]}
\newcommand{\kbrab}[1]{\big[#1\big]}

\newcommand{\absB}[1]{\Big|#1\Big|}

\begin{document}
\maketitle
\begin{abstract}
We consider a class of numerical approximations to the Caputo fractional
derivative. Our assumptions permit the use of nonuniform time steps, such as is
appropriate for accurately resolving the behavior of a solution whose
temporal derivatives are singular at~$t=0$.  The main result is a type of fractional
Gr\"{o}nwall inequality and we illustrate its use by outlining some stability and
convergence estimates of schemes for fractional reaction-subdiffusion problems.
This approach extends earlier work that used the familiar L1 approximation to
the Caputo fractional derivative, and will facilitate the analysis of higher
order and linearized fast schemes.
\end{abstract}


%%%%%%%%%%%%%%%%%%%%%%%%%%%%%%%%%%%%%%%%%%%%%%%%%%%%%%%%%%%%%%%%%%%%%%%%%%%%%%%
\section{Introduction}
This paper builds on earlier results~\cite{LiaoLiZhang:2018} for
the nonuniform L1 method applied to the time discretization of a fractional
reaction-subdiffusion problem \cite{Podlubny:1999} in a spatial domain~$\Omega$,
\begin{equation}\label{eq: reaction-diffusion}
\begin{aligned}
\fd{\alpha} u+\opL u&=f(x, t, u)&&\text{for $x\in\Omega$ and $0<t\le T$,}\\
u&=u_0(x)&&\text{for $x\in\Omega$ when $t=0$,}\\
u&=0&&\text{for $x\in\partial\Omega$ and $0<t<T$.}
\end{aligned}
\end{equation}
Here, $\fd{\alpha}={}^{\text{C}}_0\fd{\alpha}$ denotes the Caputo fractional
derivative of order~$\alpha$ with respect to time~$t$, with~$0<\alpha<1$,
and $\opL$ is a linear, second-order, strongly-elliptic partial differential
operator in the spatial variable(s)~$x$. We establish a discrete
Gr\"{o}nwall inequality intended for the error analysis of higher-order
time discretizations \cite{LiaoMcLeanZhang:2018} and linearized fast algorithms
\cite{LiaoYanZhang:2018} for solving~\eqref{eq: reaction-diffusion}
that employ nonuniform step sizes.

In any numerical methods for solving the
reaction-subdiffusion problem~\eqref{eq: reaction-diffusion}, a key consideration
is that the solution~$u(x,t)$ is typically less regular than would be the case
for a classical parabolic PDE (which arises as the limiting case
when~$\alpha\to1$).  For example, in the simplest case~$f(x,t,u)\equiv0$
when~\eqref{eq: reaction-diffusion} is linear and homogeneous, let
$\varphi_{\opL}$ be a Dirichlet eigenfunction of~$\opL$ on~$\Omega$, with
eigenvalue~$\lambda_{\opL}>0$, so that
$\opL\varphi_{\opL}=\lambda_{\opL}\varphi_{\opL}$. Let $E_\alpha$
denote the Mittag--Leffler function,
\begin{equation}\label{eq: Mittag-Leffler}
E_\alpha(z)\defeq\sum_{k=0}^\infty\frac{z^k}{\Gamma(1+k\alpha)},
\end{equation}
and choose as the initial data~$u_0(x)=\varphi_{\opL}(x)$.  Term-by-term
differentiation shows that the solution is
$u(x,t)=E_\alpha(-\lambda_{\opL} t^\alpha)\varphi_{\opL}(x)$, and so
$\partial u/\partial t=O(t^{\alpha-1})$ as~$t\to0$, whereas the solution of
the classical parabolic equation,
$u(x,t)=e^{-\lambda_{\opL}t}\varphi_{\opL}(x)$, is a smooth function of~$t$.
Sakamoto and Yamamoto~\cite{SakamotoYamamoto2011} study the (lack of) regularity
of~$u$ for more general initial data~$u_0$ and a linear source term~$f=f(x,t)$.
In fact, $u$ can only be a smooth function of~$t$ if the initial data and source
term satisfy some restrictive compatibility conditions~\cite{Stynes2016}.

The nature of polynomial interpolation means that the
convergence rate of the L1 or similar approximations to~$\fd{\alpha} u$ is
limited by the smoothness of the solution~$u$.  In the presence of a fixed
singularity at~$t=0$ of the type described above, an established technique to
restore an optimal convergence rate is to employ a graded mesh
\begin{equation}\label{eq: graded mesh}
t_n\defeq(n/N)^\gamma T\quad\text{for $0\le n\le N$,}
\end{equation}
where the parameter~$\gamma\ge1$ must be adapted to the strength of the
singularity.  Choosing $\gamma=1$ results in a uniform mesh, and the larger the
value of~$\gamma$ the more strongly the grid points are concentrated near~$t=0$.
For example, such meshes have long been used in the numerical solution of
Fredholm~\cite{Graham1982} and Volterra~\cite{Brunner1985} integral equations,
and their use for time-fractional PDEs~\cite{McLeanMustapha2007} is now well
established.

Early papers on L1 schemes \cite{LinXu2007,SunWu2006} assumed a uniform
step size~$\tau$, and showed that if $u$ is smooth then the time discretization
error is $O(\tau^{2-\alpha})$. Recently, Jin, Lazarov and
Zhou~\cite{JinLazarovZhou2016} presented a new analysis, based on generating
functions, that permitted nonsmooth initial data~$u_0$. They showed that if
$f\equiv0$~and $u_0\in L_2(\Omega)$, then the error in the norm of~$L_2(\Omega)$
due to the time discretization is $O(\tau t_n^{-1})$.  Thus, for $t_n$ bounded
away from zero, the method achieves first-order accuracy in time.
Yan, Khan and Ford \cite{YanKhanFord:2018} proposed a modified L1 scheme
and obtained error estimates for smooth and nonsmooth initial data.
It was shown that the modified L1 scheme on a uniform mesh has a convergence
rate of $O(\tau^{2-\alpha})$.
Aliknanov~\cite{Alikhanov2015} introduced the L2-$1_\sigma$ formula, a
modification of the L1 method that uses
piecewise-quadratic instead of piecewise-linear interpolation, and approximates
$\fd{\alpha} u$ at an offset grid point~$t_{j+\sigma}=(j+\sigma)\tau$.  He
showed that if $u$ is sufficiently smooth then the time discretization error is
$O(\tau^2)$ for the special choice~$\sigma=1-\alpha/2$.

Although nonuniform meshes are flexible and reasonably convenient for practical
implementation, they can significantly complicate the numerical analysis of
schemes, both with respect to stability and consistency.  Stynes, O'Riordan
and Gracia~\cite{StynesEtAl2017} considered the L1 method on a graded
mesh of the form~\eqref{eq: graded mesh} applied
to~\eqref{eq: reaction-diffusion} for the case $\opL u=-u_{xx}$ and a linear
reaction term~$f(x,t,u)=-c(x)u+g(x,t)$. They showed that, given the typical
singular behavior of~$u$, the maximum error in the fully-discrete solution is
of order~$N^{-\min\{2-\alpha,\gamma\alpha\}}$.  (Here we ignore the additional
error due to the spatial discretization.)  Thus, for a uniform mesh the
error is $O(N^{-\alpha})$, but if $\gamma=(2-\alpha)/\alpha$
then the error is $O(N^{\alpha-2})$. Their stability
analysis requires $c(x)\ge0$, which prevents extending the approach to deal with
a reaction term that is nonlinear but uniformly Lipschitz in~$u$.  This
limitation was overcome recently in the precursor~\cite{LiaoLiZhang:2018} to the
present work by exploiting a novel discrete fractional Gr\"{o}nwall inequality for the L1
method.

Nonetheless, practical applications of the discrete
Gr\"{o}nwall inequality in its basic form~\cite{LiaoLiZhang:2018} are still
limited because it does not apply to other numerical approximation schemes for
the Caputo derivative and excludes certain adaptive time meshes required
to resolve complex behaviors (physical oscillations, blowup and so on)
in nonlinear fractional differential equations.  Also, the proof relies on
specific properties of the L1 kernels~$a^{(n)}_{n-k}$ and their complementary
discrete kernels~$P^{(n)}_{n-k}$, with a key
step~\cite[Lemma~2.1]{LiaoLiZhang:2018} employing rough estimates of the
truncation error that, to a large extent, rely on the simple form
of the~$a^{(n)}_{n-k}$.  In summary, the main novel contributions of the
present work are threefold:
\begin{itemize}
\item[(i)] to generalize the discrete Gronwall inequality, permitting its use
with a variety of discretizations of the Caputo derivative, not just the L1
scheme;
\item[(ii)] to provide a concise proof based on two simple assumptions on the
discrete kernels, independent of their precise form;
\item[(iii)] to permit a more general class of nonuniform meshes or adaptive
time grids, not just the graded meshes for resolving the initial singularity.
\end{itemize}
In more detail, \cref{sec: discrete fractional} defines a discrete fractional
derivative~\eqref{eq: discrete Caputo} having the form of the classical L1
approximation but with general discrete kernels.  We formulate three
assumptions required for our theory.  The first two impose a monotonicity
property (\Ass{1}) and a lower bound (\Ass{2}) on the discrete kernels, and the
third (\Ass{3}) places a mild restriction on the local step-size ratio.  We
give some examples of schemes satisfying these assumptions, and define a family
of complementary discrete kernels, generalizing those introduced in the earlier
paper~\cite{LiaoLiZhang:2018}.  \Cref{lem: E_alpha} establishes a key
estimate involving the discrete kernels and the Mittag--Leffler
function~\eqref{eq: Mittag-Leffler}.  In \cref{sec: discrete Gronwall},
we prove our main result, a discrete fractional Gr\"{o}nwall inequality
stated as \cref{thm: gronwall}, and provide, in \cref{rem: Caputo BDF2}, a
strategy to treat cases where the monotonicity assumption breaks down.
\Cref{sec: stability consistency} illustrates the use of the Gronwall inequality
in conjunction with an abstract Galerkin method for the spatial discretization.
Finally, a short appendix proves two technical inequalities needed for the
stability analysis of \cref{sec: stability consistency}.

The generalized results proved below will allow us to show, in two
companion papers~\cite{LiaoMcLeanZhang:2018,LiaoYanZhang:2018}, that
Alikhanov's L2-1${}_\sigma$ formula can achieve second-order accuracy on certain
nonuniform time grids and that a linearized fast algorithm is unconditionally
convergent for nonlinear subdiffusion equations.

%when a graded mesh~\eqref{eq: graded mesh} is used to compensate
%for the singular behavior of~$u$.

%%%%%%%%%%%%%%%%%%%%%%%%%%%%%%%%%%%%%%%%%%%%%%%%%%%%%%%%%%%%%%%%%%%%%%%%%%%%%%%
\section{Discrete fractional derivative}\label{sec: discrete fractional}
Recall that the Riemann--Liouville fractional integral operator of
order~$\beta>0$ is defined by \cite{OldhamSpanier1974,Podlubny:1999}
\[
(\mathcal{I}^\beta v)(t)\defeq\int_0^t\omega_\beta(t-s)v(s)\zd{s}
	\quad\text{for $t>0$,}
	\quad\text{where $\omega_\beta(t)\defeq\frac{t^{\beta-1}}{\Gamma(\beta)}$,}
\]
and, in turn, the Caputo fractional derivative is defined by
\begin{equation}\label{eq: Caputo}
(\fd{\alpha}v)(t)\defeq(\mathcal{I}^{1-\alpha}v')(t)
	=\int_0^t\omega_{1-\alpha}(t-s)v'(s)\zd{s}
	\quad\text{for $t>0$.}
\end{equation}
For (possibly nonuniform) time levels
$0=t_0<t_1<t_2<\cdots<t_N=T$, we denote the $n$th step size
by~$\tau_n\defeq t_n-t_{n-1}$, fix an offset parameter~$\theta\in[0,1)$ and
define
\[
t_{n-\theta}\defeq\theta t_{n-1}+(1-\theta)t_n
\quad\text{and}\quad
v^{n-\theta}\defeq\theta v^{n-1}+(1-\theta)v^n,
\]
where~$v^k$ may be any sequence. Letting $v^k\approx v(t_k)$ and
$\diff v^k\defeq v^k-v^{k-1}$,
we consider a discrete Caputo derivative (not necessarily a direct approximation of \eqref{eq: Caputo},
see Remark \ref{rem: other Caputo})
given by a convolution-like sum, as follows,
\begin{equation}\label{eq: discrete Caputo}
(\dfd{\alpha}v)^{n-\theta}\defeq\sum_{k=1}^n
	A^{(n)}_{n-k}\diff v^k\quad\text{for $1\le n\le N$.}	
\end{equation}
Here, the corresponding discrete convolution kernels are written as~$A^{(n)}_{n-k}$ instead
of~$A_{nk}$ to reflect the convolution structure of the fractional derivative.
%In some cases, a special choice of~$\theta$ results in higher accuracy~\cite{Alikhanov2015}.
Our theory requires the following three
assumptions:
\begin{description}
\item[A1.] The discrete kernels are positive and monotone, that is,
\[
A^{(n)}_0\ge A^{(n)}_1\ge A^{(n)}_2\ge\cdots A^{(n)}_{n-1}>0
\quad\text{for $1\le n\le N$.}
\]
\item[A2.] There is a constant~$\pi_A>0$ such that the discrete kernels satisfy
the lower bound
\[
A^{(n)}_{n-k}\ge\frac{1}{\pi_A\tau_k}\int_{t_{k-1}}^{t_k}
\omega_{1-\alpha}(t_n-s)\zd{s}\quad\text{for $1\le k\le n\le N$.}
\]
\item[A3.] There is a constant~$\rhomax>0$ such that the step size
ratios~$\rho_k\defeq\tau_k/\tau_{k+1}$ satisfy
\[
\rho_k\le\rhomax\quad\text{for $1\le k\le N-1$.}
\]
\end{description}

The boundedness and monotonicity assumptions \Ass{1} and \Ass{2}
on the discrete convolution kernels~$A^{(n)}_{n-k}$ are valid for several
frequently-used discrete Caputo derivatives, at least if assumption~\Ass{3} is
satisfied for appropriate~$\rho$.  Included are the well-known L1 formula
\cite{LiaoLiZhang:2018,LinXu2007,OldhamSpanier1974,StynesEtAl2017,SunWu2006}, the
fast L1 formula~\cite{LiaoYanZhang:2018},
the Alikhanov approximation~\cite{Alikhanov2015,LiaoZhaoTeng2016,LiaoMcLeanZhang:2018},
and their applications for multi-term and distributed-order Caputo derivatives
(see Remark \ref{rem: other Caputo}). Here we list three examples on nonuniform grids.
Note that, the local mesh parameter~$\rho$ from~\Ass{3} will always
appear in our discrete fractional Gr\"{o}nwall inequality and our stability estimates.

\begin{example}[nonuniform L1 formula]\label{exam: L1}
The widespread L1 formula~\cite[p.~140]{OldhamSpanier1974}
uses $\theta=0$~and $v'(s)\approx\diff v^k/\tau_k$ (linear interpolation) to obtain
\[
(\dfd{\alpha}v)^n\defeq\sum_{k=1}^na^{(n)}_{n-k}\diff v^k\quad \text{with} \quad
a^{(n)}_{n-k}\defeq\frac{1}{\tau_k}\int_{t_{k-1}}^{t_k}
    \omega_{1-\alpha}(t_n-s)\zd{s}.
\]
This sum has the desired form~\eqref{eq: discrete Caputo} where (using the
integral mean value theorem),
\begin{equation}\label{eq: L1 weights}
A^{(n)}_{n-k}:=a^{(n)}_{n-k}=\omega_{1-\alpha}(t_n-s_{nk})\quad \text{for some $s_{nk}\in[t_{k-1},t_k]$.}
\end{equation}
It follows that assumption~\Ass{1} is satisfied, and
\Ass{2} holds with~$\pi_A=1$.
\end{example}
\begin{example}[fast L1 formula]\label{exam: fast L1}
In the two-level fast L1 approximation \cite{LiaoYanZhang:2018},
the sum-of-exponentials technique is applied to approximate the weakly singular
kernel $\omega_{1-\alpha}(t-s)$. That is, for a user-given absolute tolerance
error $\epsilon\ll1$ and a cut-off time $\Delta{t}>0$, one determines a
positive integer $N_{q}$, positive quadrature nodes $\theta^{\ell}$ and positive
weights $\varpi^{\ell}$ $(1\leq \ell\leq N_q)$ such that
\[
\absB{\omega_{1-\alpha}(t_k-s)-\sum_{\ell=1}^{N_q}\varpi^{\ell}e^{-\theta^{\ell}(t_k-s)}}
\leq \epsilon \quad \forall\, t_k\in[s+\Delta{t},T].
\]
Then we use $\theta=0$~and $v'(s)\approx\diff v^k/\tau_k$ (linear interpolation) to obtain
\[
 \brat{\ffd{\alpha}u}^n\defeq a_{0}^{(n)}\nabla_{\tau}u^n
  +\sum_{\ell=1}^{N_q}\varpi^{\ell}e^{-\theta^{\ell}\tau_n}H^{\ell}(t_{n-1}),
\quad n\geq1,
\]
where $H^{\ell}(t_{k})$ satisfies $H^{\ell}(t_{0})=0$ and the
recurrence relationship
\[
H^{\ell}(t_{k})=e^{-\theta^{\ell}\tau_{k}}H^{\ell}(t_{k-1})
    +\frac{1}{\tau_k}\int_{t_{k-1}}^{t_{k}}
    e^{-\theta^{\ell}(t_{k}-s)}\nabla_{\tau}u^{k}\zd s,\quad
    k\geq1, \;1\leq \ell\leq N_q\,.
\]
This approximation also has the form~\eqref{eq: discrete Caputo} with $\theta=0$,
\begin{align*}
 A_{0}^{(n)}\defeq a_{0}^{(n)} \quad\text{and}\quad
A_{n-k}^{(n)}\defeq\frac{1}{\tau_k}\int_{t_{k-1}}^{t_{k}}\sum_{\ell=1}^{N_q}\varpi^{\ell}e^{-\theta^{\ell}(t_{n}-s)}\zd s
\quad \text{for $1\le k\leq n-1$}.
\end{align*}
If the tolerance error $\epsilon$ is small enough such that
$\epsilon\leq\min\left\{\frac{1}{3}\omega_{1-\alpha}(T),\alpha\,\omega_{2-\alpha}(1)\right\}$,
then
\cite[Lemma 2.5]{LiaoYanZhang:2018} ensures that \Ass{1}-\Ass{2} hold true with $\pi_A=3/2$.
\end{example}
\begin{example}[nonuniform Alikhanov formula]\label{exam: Alikhanov}
Let $\Pi_{1,k}v$ be the linear interpolant of a function~$v$ with
respect to the nodes $t_{k-1}$~and $t_k$, and let $\Pi_{2,k}v$ denote the
quadratic interpolant with respect to $t_{k-1}$, $t_k$~and $t_{k+1}$.
Taking a special choice $\theta=\alpha/2$, and applying the linear and
quadratic polynomial interpolations, we have the nonuniform Alikhanov formula
~\cite{LiaoZhaoTeng2016,LiaoMcLeanZhang:2018}
\begin{align*}
(\dfd{\alpha}v)^{n-\theta}
	\defeq&\,\sum_{k=1}^{n-1}\int_{t_{k-1}}^{t_k}\omega_{1-\alpha}(t_{n-\theta}-s)\bra{\Pi_{2,k}v}'(s)\zd{s}\\
&\,+\int_{t_{n-1}}^{t_{n-\theta}}\omega_{1-\alpha}(t_{n-\theta}-s)
		\bra{\Pi_{1,n}v}'(s)\zd{s}\quad\text{for $n\ge1$}.
\end{align*}
This formula can be written as the form~\eqref{eq: discrete Caputo} with
$A_0^{(1)}\defeq \hat{a}_0^{(1)}$ for $n=1$ and, for $n\geq2$,
\begin{equation*}
A^{(n)}_{n-k}\defeq\begin{cases}
	\hat{a}^{(n)}_0+\rho_{n-1}\hat{b}^{(n)}_1,
	&\text{for $k=n$,}\\
	\hat{a}^{(n)}_{n-k}+\rho_{k-1}\hat{b}^{(n)}_{n-k+1}-\hat{b}^{(n)}_{n-k},
	&\text{for $2\le k\le n-1$,}\\
	\hat{a}^{(n)}_{n-1}-\hat{b}^{(n)}_{n-1},
	&\text{for $k=1$,}
\end{cases}
\end{equation*}
where the discrete coefficients
$\hat{a}_{n-k}^{(n)}$ and $\hat{b}_{n-k}^{(n)}$ are defined by
\begin{align*}
&\hat{a}^{(n)}_0\defeq\frac1{\tau_n}\int_{t_{n-1}}^{t_{n-\theta}}
	\omega_{1-\alpha}(t_{n-\theta}-s)\zd{s},\\
&\hat{a}^{(n)}_{n-k}\defeq\frac{1}{\tau_k}\int_{t_{k-1}}^{t_k}\omega_{1-\alpha}(t_{n-\theta}-s)\zd{s}\quad \text{for $1\le k\leq n-1$},\\
&\hat{b}^{(n)}_{n-k}\defeq\frac{2}{\tau_k(\tau_k+\tau_{k+1})}\int_{t_{k-1}}^{t_k}
	(s-t_{k-\frac12})\omega_{1-\alpha}(t_{n-\theta}-s)\zd{s}\quad \text{for $1\le k\leq n-1$}.
\end{align*}
The theoretical properties in \cite[Theorem 2.2]{LiaoMcLeanZhang:2018} assure
\Ass{1}--\Ass{2} with $\pi_A=11/4$ provided the local mesh assumption \Ass{3}
holds with the maximum step size ratio $\rho=7/4$.
\end{example}

%Typically, the third assumption \Ass{3} is satisfied by the graded mesh~\eqref{eq: graded mesh} with~$\rhomax=1$.

We now continue to introduce an important tool: the complementary discrete convolution kernels.
The semigroup property of the fractional integral,
$\mathcal{I}^\alpha\mathcal{I}^\beta=\mathcal{I}^{\alpha+\beta}$, holds because
the integral kernels satisfy $\omega_\alpha*\omega_\beta=\omega_{\alpha+\beta}$.  It
follows that
\begin{equation}\label{eq: omega identity}
\int_s^t\omega_\alpha(t-\mu)\omega_{1-\alpha}(\mu-s)\zd\mu=\omega_1(t-s)=1
	\quad\text{for all $0<s<t<\infty$,}
\end{equation}
and it motives us to seek a family
of complementary discrete convolution kernels~$P^{(n)}_{n-j}$ having the identical property
\begin{equation}\label{eq: P weights}
\sum_{j=m}^n P^{(n)}_{n-j}A^{(j)}_{j-m}\equiv1
	\quad\text{for $1\le m\le n\le N$.}
\end{equation}
In fact, by taking $m=k$ and $m=k+1$,
\[
P^{(n)}_{n-k}A^{(k)}_0+\sum_{j=k+1}^nP^{(n)}_{n-j}A^{(j)}_{j-k}
	=1=\sum_{j=k+1}^nP^{(n)}_{n-j}A^{(j)}_{j-(k+1)} \;,\quad 1\le k\le n-1,
\]
we see that
\[
P^{(n)}_{n-k}=\frac{1}{A^{(k)}_0}\sum_{j=k+1}^n P^{(n)}_{n-j}\braB{
	A^{(j)}_{j-k-1}-A^{(j)}_{j-k}},\quad 1\le k\le n-1,
\]
and the complementary discrete kernels may be defined via the
recursion~\cite{LiaoLiZhang:2018}
\begin{equation}\label{eq: P recursion}
P^{(n)}_0\defeq\frac{1}{A^{(n)}_0},\;
P^{(n)}_j\defeq\frac{1}{A^{(n-j)}_0}\sum_{k=0}^{j-1}\braB{
	A^{(n-k)}_{j-k-1}-A^{(n-k)}_{j-k}}P^{(n)}_k
	\;\;\text{for $1\le j\le n-1$.}
\end{equation}

\begin{example}[Pictures of $A_j^{(n)}$ and $P_j^{(n)}$ of L1 formula]\label{exam: L1 figures}
Consider the widespread L1 approximation in Example \ref{exam: L1}.
Figure~\ref{fig: weights} plots the L1 discrete kernels $A^{(n)}_{j}$ and the complementary discrete kernels~$P^{(n)}_{j}$
 when $T=1$~and $n=30$ for three graded meshes of the form~\eqref{eq: graded mesh}.
 \begin{figure}
\caption{Top: the L1 discrete kernel~\eqref{eq: L1 weights} for three
different meshes of the form~\eqref{eq: graded mesh} in the case $T=1$~and
$n=30$. Bottom: the complementary discrete kernels~$P^{(n)}_j$.}
\label{fig: weights}
\begin{center}
\includegraphics[scale=0.50]{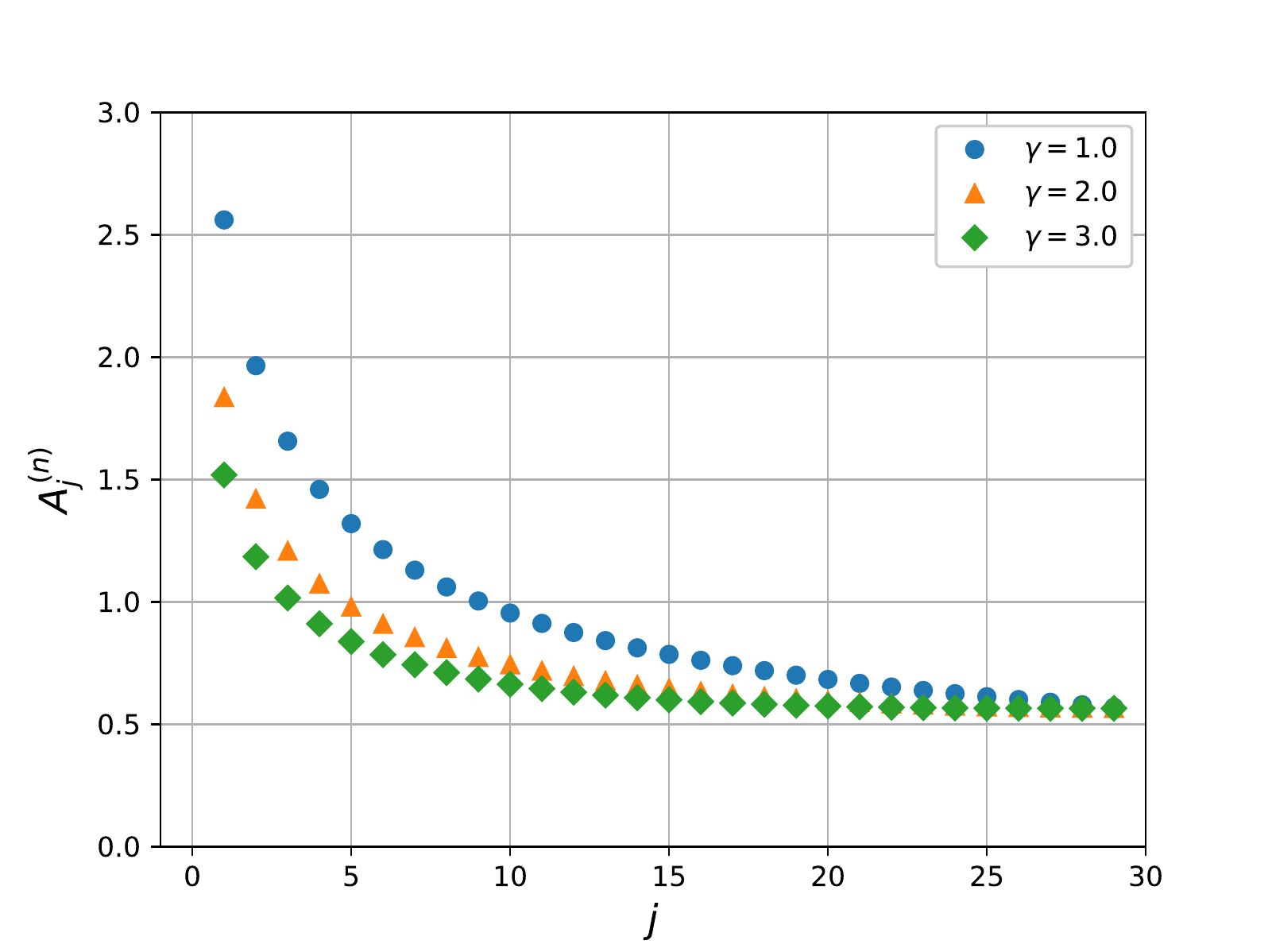}
\includegraphics[scale=0.50]{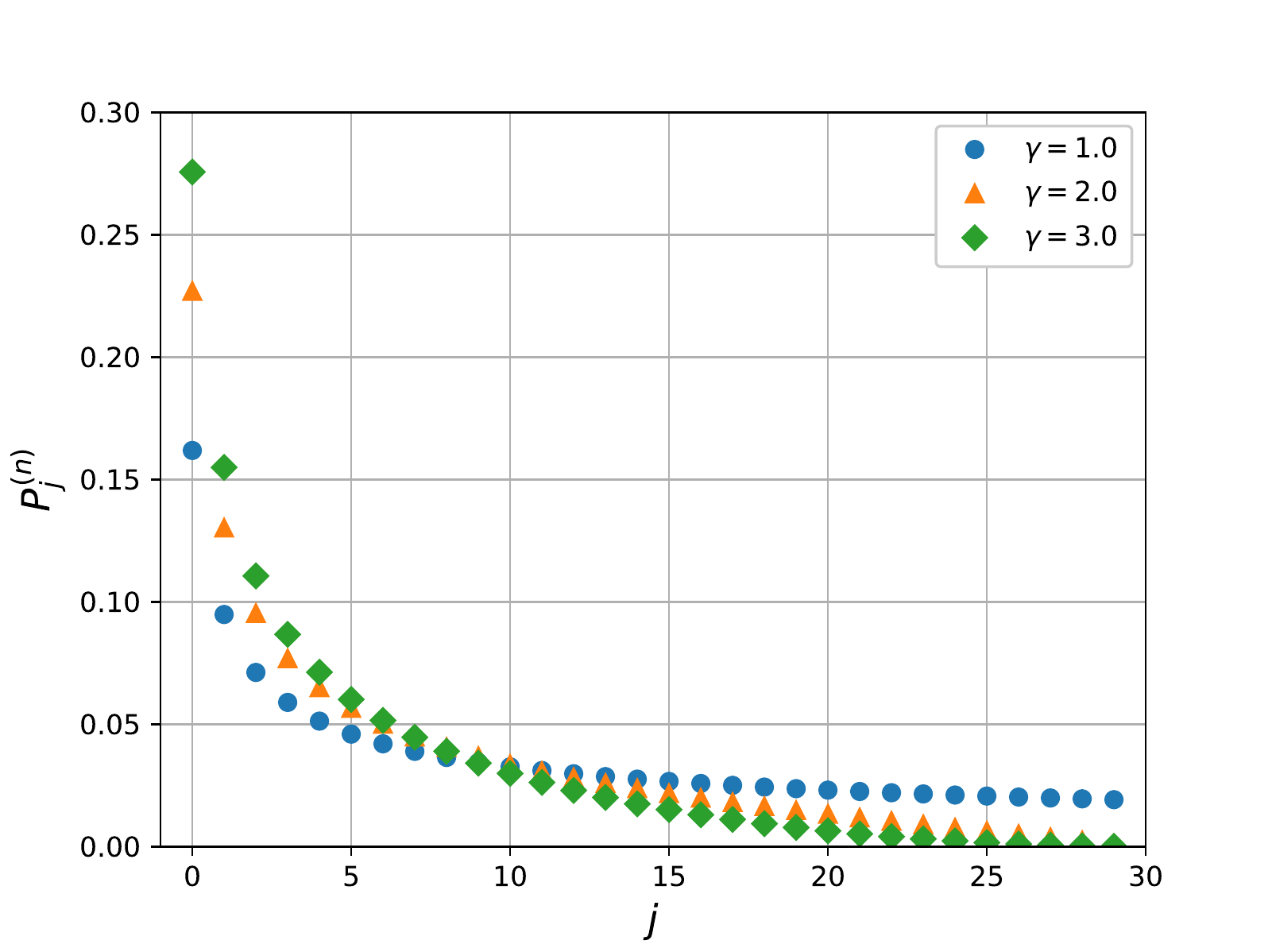}
\end{center}
\end{figure}
\end{example}

As a consequence of the identity~\eqref{eq: omega identity}, we find that
\[
\int_0^t\omega_\alpha(t-s)(\fd{\alpha} v)(s)\zd{s}=\int_0^t v'(s)\zd{s},
\]
which provides the inspiration for the second part of the next lemma.

\begin{lemma}\label{lem: Pn A1 A2}
Let the assumptions \Ass{1}~and \Ass{2} hold.
\begin{enumerate}
\item The discrete kernels $P^{(n)}_{j}$
in~\eqref{eq: P recursion} having the property \eqref{eq: P weights} satisfy
\[
0\le P^{(n)}_{n-j}\le\pi_A\Gamma(2-\alpha)\tau_n^\alpha
	\quad\text{for $1\le j\le n\le N$,}
\]
and
\begin{align}
&\sum_{j=1}^{n}P_{n-j}^{(n)}\,\omega_{1-\alpha}(t_j)\leq \pi_A\quad\text{for $1\le n\le N$.}\label{eq: sum P1}
%\\&\sum_{j=1}^nP^{(n)}_{n-j}\le\pi_A\omega_{1+\alpha}(t_n)
%    \quad\text{for $1\le n\le N$.}\label{eq: sum P}
\end{align}
\item If $v:[0,T]\to\R$ is any continuous, piecewise-$C^1$ function such that
$v'$ is non-negative and monotone decreasing, then
\begin{align} \label{lem22}
\sum_{j=1}^n P^{(n)}_{n-j}(\fd{\alpha} v)(t_j)\le\pi_A\int_0^{t_n}v'(s)\zd{s} \quad\text{for $1\le n\le N$.}
\end{align}
\end{enumerate}
\end{lemma}

\begin{proof}
It follows at once from  the monotonicity assumption~\Ass{1} that
$A^{(n)}_0>0$~and $A^{(n-k)}_{j-k-1}-A^{(n-k)}_{j-k}\ge0$ for $0\le k\le j-1$.
The lower bound $P^{(n)}_j\ge0$ is then clear from the
recursion~\eqref{eq: P recursion}.  Since all the discrete kernels are
non-negative, we have
\[
P^{(n)}_{n-k}A^{(k)}_0\le\sum_{j=k}^nP^{(n)}_{n-j}A^{(j)}_{j-k}=1
\]
and taking $n=k$ in the assumption~\Ass{2} gives
\[
A^{(k)}_0\ge\frac{1}{\pi_A\tau_k}\int_{t_{k-1}}^{t_k}
	\omega_{1-\alpha}(t_k-s)\zd{s}
	=\frac{\omega_{2-\alpha}(\tau_k)}{\pi_A\tau_k}
	=\frac{1}{\Gamma(2-\alpha)\pi_A\tau_k^\alpha},
\]
so the complementary discrete convolution kernels $P^{(n)}_{n-k}$ are well-defined and satisfy
the upper bound
$P^{(n)}_{n-k}\le1/A^{(k)}_0\le\Gamma(2-\alpha)\pi_A\tau_k^\alpha$.
Furthermore, the assumption {\bf A2} and the identity \eqref{eq: P weights} imply that
$\omega_{1-\alpha} (t_j) \leq \pi_A A_{j-1}^{(j)}$ and
\begin{align*}
&\sum_{j=1}^{n}P_{n-j}^{(n)}\,\omega_{1-\alpha}(t_j)\leq \pi_A \sum_{j=1}^{n}P_{n-j}^{(n)} A_{j-1}^{(j)} = \pi_A\quad \text{for $n\geq1$},
\end{align*}
which completes the proof of part~1.

Recall Chebyshev's sorting
inequality~\cite[p.~168, item~236.]{HardyLittlewoodPolya1934}: if $f$ is
monotone increasing and $g$ is monotone decreasing on the interval~$[a,b]$, and
if both functions are integrable, then
\[
(b-a)\int_a^bf(s)g(s)\zd{s}\le\int_a^b f(t)\zd t\int_a^b g(s)\zd{s}.
\]
Taking $[a,b]=[t_{k-1},t_k]$, $f(s)=\omega_{1-\alpha}(t_j-s)$~and $g(s)=v'(s)\ge0$,
and using \Ass{2}, we see that
\begin{align*}
(\fd{\alpha} v)(t_j)&=\sum_{k=1}^j\int_{t_{k-1}}^{t_k}
	\omega_{1-\alpha}(t_j-s)v'(s)\zd{s}\\
	&\le\sum_{k=1}^j\frac{1}{\tau_k}
		\int_{t_{k-1}}^{t_k}\omega_{1-\alpha}(t_j-t)\zd t
		\int_{t_{k-1}}^{t_k}v'(s)\zd{s}
	\le\pi_A\sum_{k=1}^jA^{(j)}_{j-k}
		\int_{t_{k-1}}^{t_k}v'(s)\zd{s}.
\end{align*}
Thus, from the identical property~\eqref{eq: P weights} of the discrete kernels~$P^{(n)}_{n-j}$, we conclude that
\begin{align*}
\sum_{j=1}^nP^{(n)}_{n-j}(\fd{\alpha} v)(t_j)
	&\le\sum_{j=1}^nP^{(n)}_{n-j}\pi_A\sum_{k=1}^jA^{(j)}_{j-k}
		\int_{t_{k-1}}^{t_k}v'(s)\zd{s}\\
	&=\pi_A\sum_{k=1}^n\int_{t_{k-1}}^{t_k}v'(s)\zd{s}
	\sum_{j=k}^nP^{(n)}_{n-j}A^{(j)}_{j-k}
	=\pi_A\sum_{k=1}^n\int_{t_{k-1}}^{t_k}v'(s)\zd{s},
\end{align*}
and part~2 follows.  %Finally, by taking $v(t)=\omega_{1+\alpha}(t)$,
%so that $\fd{\alpha} v(t)=1$, we obtain
%\[
%\sum_{j=1}^nP^{(n)}_{n-j}=
%    \sum_{j=1}^nP^{(n)}_{n-j}\fd{\alpha}\omega_{1+\alpha}(t_j)
%    \le\pi_A\int_0^{t_n}\omega_{1+\alpha}'(s)\zd{s}
%    =\pi_A\omega_{1+\alpha}(t_n),
%\]
\end{proof}

When~\Ass{3} also holds, we have a variant of the second part of
\cref{lem: Pn A1 A2}.

\begin{lemma}\label{lem: A1-A3}
Let the assumptions \Ass{1}--\Ass{3} hold.  If $v:[0,T]\to\R$ is any continuous,
piecewise-$C^1$ function
such that $v'$ is non-negative and monotone, then
\[
\sum_{j=1}^{n-1}P^{(n)}_{n-j}(\fd{\alpha} v)(t_j)
	\le\max(1,\rhomax)\pi_A\int_0^{t_n}v'(s)\zd{s}\quad\text{for $1\le n\le N$.}
\]
\end{lemma}
\begin{proof}
If $v'$ is non-negative and monotone \emph{de}creasing, then $\fd{\alpha} v(t_j)\ge0$ and the results of
\cref{lem: Pn A1 A2} imply that
\[
\sum_{j=1}^{n-1}P^{(n)}_{n-j}(\fd{\alpha} v)(t_j)
	\le\sum_{j=1}^nP^{(n)}_{n-j}(\fd{\alpha}v)(t_j)
	\le\pi_A\int_0^{t_n}v'(s)\zd{s}\,.
\]
 Otherwise, if $v'$ is
monotonely \emph{increasing}, then
\begin{align*}
\sum_{j=1}^{n-1}P^{(n)}_{n-j}(\fd{\alpha} v)(t_j)
	&=\sum_{j=1}^{n-1}P^{(n)}_{n-j}\sum_{k=1}^j\int_{t_{k-1}}^{t_k}
		\omega_{1-\alpha}(t_j-s)v'(s)\zd{s}\\
	&\le\sum_{j=1}^{n-1}P^{(n)}_{n-j}\sum_{k=1}^j v'(t_k)
	\int_{t_{k-1}}^{t_k}\omega_{1-\alpha}(t_j-s)\zd{s}\\
	&\le\pi_A\sum_{j=1}^{n-1}P^{(n)}_{n-j}
		\sum_{k=1}^j v'(t_k)\tau_kA^{(j)}_{j-k}\\
	&=\pi_A\sum_{k=1}^{n-1} v'(t_k)\tau_k
		\sum_{j=k}^{n-1}P^{(n)}_{n-j}A^{(j)}_{j-k}
	\le\pi_A\sum_{k=1}^{n-1}v'(t_k)\tau_k\\
	&\le\pi_A\sum_{k=1}^{n-1}v'(t_k)\rho_k\tau_{k+1}
	\le\rhomax\pi_A\sum_{k=1}^{n-1}\int_{t_k}^{t_{k+1}}v'(s)\zd{s},
\end{align*}
and the desired estimate again holds.
\end{proof}

We can use \cref{lem: A1-A3} to prove the following property of the
Mittag--Leffler function~\eqref{eq: Mittag-Leffler}.

\begin{lemma}\label{lem: E_alpha}
Let the assumptions \Ass{1}--\Ass{3} hold. For any real $\mu>0$,
\[
\sum_{j=1}^{n-1}P^{(n)}_{n-j}E_\alpha(\mu t_j^\alpha)
	\le\pi_A\max(1,\rhomax)\frac{E_\alpha(\mu t_n^\alpha)-1}{\mu}
	\quad\text{for $1\le n\le N$.}
\]
\end{lemma}
\begin{proof}
The series definition~\eqref{eq: Mittag-Leffler} shows that
\[
E_\alpha(\mu t^\alpha)
	=1+\sum_{k=1}^\infty\frac{\mu^kt^{k\alpha}}{\Gamma(1+k\alpha)}
        =1+\sum_{k=1}^\infty\mu^kv_k(t),
\]
where $v_k(t)=\omega_{1+k\alpha}(t)$ and we have $v_k'(t)=\omega_{k\alpha}(t)>0$ for all~$k\ge1$.
If $1\le k\le1/\alpha$, then $-1\le k\alpha-1\le0$ and
$v_k''(t)=\omega_{k\alpha-1}(t)\le0$ for all~$t>0$.  Otherwise, if $k>1/\alpha$,
then $k\alpha-1>0$ and $v_k''(t)>0$ for all~$t>0$.  Thus, $v_k'$ is always
non-negative and monotone, so we may apply \cref{lem: A1-A3} and deduce
that
\[
\sum_{j=1}^{n-1}P^{(n)}_{n-j}(\fd{\alpha} v_k)(t_j)
    \le\max(1,\rhomax)\pi_A\int_0^{t_n}v_k'(s)\zd{s}
	=\max(1,\rhomax)\pi_Av_k(t_n)\quad\text{for $k\ge1$.}
\]
Multiplying both sides of this inequality by~$\mu^k$, summing over the index~$k$, and using
the fact that
\[
\fd{\alpha} v_k(t)=\int_0^t\omega_{1-\alpha}(t-s)\omega_{k\alpha}(s)\zd{s}
	=\omega_{1+(k-1)\alpha}(t)=v_{k-1}(t)\quad\text{for all $k\ge1$,}
\]
we have
\[
\sum_{k=1}^m\mu^k\sum_{j=1}^{n-1}P^{(n)}_{n-j}v_{k-1}(t)
	\le\max(1,\rhomax)\pi_A\sum_{k=1}^m\mu^kv_k(t_n).
\]
Because the series $\sum_{k=1}^\infty\mu^kv_k(t)$ is absolutely convergent and $\omega_1(t)=1$,
the desired inequality follows after interchanging the
sums on the left-hand side and then sending $m\to\infty$. The proof is completed.
\end{proof}

%%%%%%%%%%%%%%%%%%%%%%%%%%%%%%%%%%%%%%%%%%%%%%%%%%%%%%%%%%%%%%%%%%%%%%%%%%%%%%%
\section{Discrete fractional Gr\"{o}nwall inequality}\label{sec: discrete Gronwall}
Our main result is stated in the next theorem.
The proof is similar to that of \cite[Lemma 2.2]{LiaoLiZhang:2018},
but we include it here to incorporate the nonuniform mesh parameter $\rho$ in
\Ass{3}, which does not appear in discrete Gr\"{o}nwall inequalities
for classical parabolic equations.
\begin{theorem}\label{thm: gronwall}
Let the assumptions \Ass{1}--\Ass{3} hold, let $0\le\theta<1$, and
let $(g^n)_{n=1}^N$ and $(\lambda_l)_{l=0}^{N-1}$ be given non-negative
sequences.  Assume further that there exists a constant $\Lambda$
(independent of the step sizes) such that
$\Lambda\ge\sum_{l=0}^{N-1}\lambda_l$, and that the maximum step size satisfies
\[
\max_{1\le n\le N}\tau_n
    \le\frac{1}{\sqrt[\alpha]{2\pi_A\Gamma(2-\alpha)\Lambda}}\,.
\]
Then, for any non-negative sequence~$(v^k)_{k=0}^N$ such that
\begin{equation}\label{eq: first Gronwall}
\sum_{k=1}^nA^{(n)}_{n-k}\diff\brab{v^k}^2\le
	\sum_{k=1}^n\lambda_{n-k}\brab{v^{k-\theta}}^2+v^{n-\theta}g^n
	\quad\text{for $1\le n\le N$,}
\end{equation}
it holds that
\begin{equation}\label{eq: gronwall conclusion}
v^n\le2E_\alpha\brab{2\max(1,\rhomax)\pi_A \Lambda t_n^\alpha}
	\biggl(v^0+\max_{1\le k\le n}\sum_{j=1}^k P^{(k)}_{k-j}g^j
	\biggr)\quad\text{for $1\le n\le N$.}
\end{equation}
\end{theorem}
\begin{proof}
We replace the index~$n$ with~$j$ in~\eqref{eq: first Gronwall}, then multiply
by~$P^{(n)}_{n-j}$ and sum over~$j$ to obtain
\begin{equation}\label{eq: gronwall 1}
\sum_{j=1}^nP^{(n)}_{n-j}\sum_{k=1}^jA^{(j)}_{j-k}\diff\brab{v^k}^2
	\le
	\sum_{j=1}^nP^{(n)}_{n-j}\sum_{k=1}^j\lambda_{j-k}\brab{v^{k-\theta}}^2+\sum_{j=1}^nP^{(n)}_{n-j}v^{j-\theta}g^j.
\end{equation}
On the left-hand side, we exchange the order of summation and use the
identity~\eqref{eq: P weights} to get
\begin{equation}\label{eq: gronwall 2}
\begin{aligned}
\sum_{j=1}^nP^{(n)}_{n-j}\sum_{k=1}^jA^{(j)}_{j-k}\diff\brab{v^k}^2
	&=\sum_{k=1}^n\diff\brab{v^k}^2\sum_{j=k}^nP^{(n)}_{n-j}A^{(j)}_{j-k}\\
	&=\sum_{k=1}^n\diff\brab{v^k}^2=(v^n)^2-(v^0)^2.
\end{aligned}
\end{equation}
Thus, it follows from \eqref{eq: gronwall 1} that
\begin{equation}\label{eq: vn2}
\brab{v^n}^2\le\brab{v^0}^2+\sum_{j=1}^nP^{(n)}_{n-j}\sum_{k=1}^j\lambda_{j-k}\brab{v^{k-\theta}}^2
+\sum_{j=1}^nP^{(n)}_{n-j}v^{j-\theta}g^j	,
\end{equation}

For brevity, let us write the claimed estimate \eqref{eq: gronwall conclusion} as $v^n\le F_nG_n$ where
\[
F_n\defeq 2E_\alpha\bigl(2\max(1,\rhomax)\pi_A\Lambda t_n^\alpha\bigr)
\quad\text{and}\quad
G_n\defeq v^0+\max_{1\le k\le n}\sum_{j=1}^k P^{(k)}_{k-j}g^j.
\]
We will use complete induction, noting that the Mittag--Leffler function \eqref{eq: Mittag-Leffler}
satisfies $E_\alpha(0)=1$ and $E_\alpha'(z)>0$ for all real~$z>0$, so $F_{n}\geq F_{n-1}\ge2$ for $n\ge2$.

If $v^1\le v^0$, then $v^1\le G_1\le F_1G_1$, as required.
Otherwise, if $v^1>v^0$, then $v^{1-\theta}\le v^1$. One deduces from \eqref{eq: vn2} that
\begin{align*}
\brab{v^{1}}^2&\le\brab{v^{0}}^2+P^{(1)}_0 v^{1-\theta}g^1
    +P^{(1)}_0\lambda_0\brab{v^{1-\theta}}^2\\
	&\le v^1\brab{v^0+P^{(1)}_0g^1}
		+P^{(1)}_0\lambda_0\brab{v^{1}}^2=v^1G_1+P^{(1)}_0\lambda_0\brab{v^{1}}^2.
\end{align*}
Part~1 of \cref{lem: Pn A1 A2} and the given restriction on the
maximum time-step imply that
\begin{equation}\label{eq: P lambda}
P^{(1)}_0\lambda_0\le\pi_A\Gamma(2-\alpha)\tau_1^\alpha\Lambda\le1/2.
\end{equation}
Thus, $\brat{v^{1}}^2\le2v^1G_1$~and so $v^1\le2G_1\le F_1G_1$,
which implies that the desired estimate holds for~$n=1$.

For the inductive step, let $2\le n\le N$ and assume that
\begin{equation}\label{eq: hypothesis}
v^k\le F_kG_k
	\quad\text{for $1\le k\le n-1$.}
\end{equation}
Choose some $k(n)$ such that $v^{k(n)}=\max_{0\le j\le n-1}v^j$.  If
$v^n\le v^{k(n)}$ then, since $F_k$~and $G_k$ are monotone increasing in~$k$,
\[
v^n\le v^{k(n)}\le F_{k(n)}G_{k(n)}\le F_nG_n,
\]
as required.  Otherwise, if $v^n>v^{k(n)}$, then
$v^{j-\theta}\le\max(v^{j-1},v^j)\le v^n$ for~$1\le j\le n$. We deduce
from~\eqref{eq: vn2} that
\begin{equation}\label{eq: gronwall 3}
\brab{v^n}^2\le v^nv^0+v^n\sum_{j=1}^n P^{(n)}_{n-j}g^j		
+v^n\sum_{j=1}^{n-1}P^{(n)}_{n-j}\sum_{k=1}^j\lambda_{j-k}v^{k-\theta}
+\brab{v^n}^2P^{(n)}_0\sum_{k=1}^n\lambda_{n-k}.
\end{equation}
Using part~1 of \cref{lem: Pn A1 A2},
\begin{equation}\label{eq: gronwall 4}
P^{(n)}_0\sum_{k=1}^n\lambda_{n-k}
\le\pi_A\Gamma(2-\alpha)\Lambda\tau_n^\alpha,
\end{equation}
so the limitation on the maximum step size implies that
\begin{equation}\label{eq: gronwall 5}
\bra{v^n}^2\le v^n\braB{G_n
+\sum_{j=1}^{n-1}P^{(n)}_{n-j}\sum_{k=1}^j\lambda_{j-k}v^{k-\theta}}
	+\frac{1}{2}\bra{v^n}^2.
\end{equation}
Thus, applying the induction hypothesis~\eqref{eq: hypothesis}, we deduce
from~\eqref{eq: gronwall 5} that
\begin{align*}
v^n&\le2G_n+2\sum_{j=1}^{n-1}P^{(n)}_{n-j}\sum_{k=1}^j\lambda_{j-k}
	\kbrab{\theta v^{k-1}+(1-\theta)v^k}\\
	&\le2G_n+2\sum_{j=1}^{n-1}P^{(n)}_{n-j}\sum_{k=1}^j\lambda_{j-k}
	\kbrab{\theta F_{k-1}G_{k-1}+(1-\theta)F_kG_k}\\
	&\le2G_n+2\sum_{j=1}^{n-1}P^{(n)}_{n-j}\sum_{k=1}^j\lambda_{j-k}F_kG_k
	\le2G_n+2\sum_{j=1}^{n-1}P^{(n)}_{n-j}F_jG_j\sum_{k=1}^j\lambda_{j-k}\\
	&\le2G_n+4\Lambda G_{n-1}\sum_{j=1}^{n-1}P^{(n)}_{n-j}
		E_\alpha\bigl(2\max(1,\rhomax)\pi_A\Lambda t_j^\alpha\bigr).
\end{align*}
Finally, by~\cref{lem: E_alpha} with~$\mu=2\max(1,\rhomax)\pi_A\Lambda$,
\[
v^n\le 2G_n+2\max(1,\rhomax)\pi_A\Lambda G_n\,
	\frac{E_\alpha\brab{2\max(1,\rhomax)\pi_A\Lambda t_n^\alpha}-1}{\max(1,\rhomax)\pi_A\Lambda}
	= F_nG_n,
\]
which completes the inductive step and the proof.
\end{proof}

\begin{remark}\label{rem: sum P}
One may use the inequality \eqref{eq: sum P1} in part 1 of
\cref{lem: Pn A1 A2} to bound the convolutional summation
$\sum_{j=1}^k P^{(k)}_{k-j}g^j$, that is,
\begin{equation*}
\sum_{j=1}^k P^{(k)}_{k-j}g^j\leq \sum_{j=1}^k P^{(k)}_{k-j}\omega_{1-\alpha}(t_j)\max_{1\leq j\leq k}\frac{g^j}{\omega_{1-\alpha}(t_j)}
\leq \pi_A\max_{1\leq j\leq k}\frac{g^j}{\omega_{1-\alpha}(t_j)}
\end{equation*}
So the discrete solution of \eqref{eq: first Gronwall} can also be bounded by
\[
v^n\le2E_\alpha\brab{2\max(1,\rhomax)\pi_A \Lambda t_n^\alpha}
	\braB{v^0+\pi_A\Gamma(1-\alpha)\max_{1\le j\le n}\{t_j^{\alpha}g^j\}}\quad\text{for $1\le n\le N$.}
\]
On the other hand, if the given sequence $(\lambda_l)_{l=0}^{N-1}$ is non-positive and the constant $\Lambda\le0$,
a similar argument will show that the discrete inequality \eqref{eq: gronwall conclusion} holds in a simpler form,
requiring only the assumptions \Ass{1}-\Ass{2} but no restrictions on time steps,
\begin{equation}\label{eq: gronwall conclusion simple}
v^n\le v^0+\max_{1\le k\le n}\sum_{j=1}^k P^{(k)}_{k-j}g^j\le v^0+\pi_A\Gamma(1-\alpha)\max_{1\le j\le n}\{t_j^{\alpha}g^j\}
	\quad\text{for $1\le n\le N$.}
\end{equation}
\end{remark}

\begin{remark}\label{rem: lambda setting}
By including the non-negative sequence $(\lambda_l)_{l=0}^{N-1}$ in
\eqref{eq: first Gronwall}, we are able to treat various numerical approaches
to solving linear and nonlinear subdiffusion problems.
Typically, the sequence takes only a few nonzero values.  Recent examples
include $\lambda_l=0$ for~$l\geq1$ in the time-weighted method from
\cref{sec: stability consistency},
and $\lambda_l=0$ for~$l\geq2$ in the one-step linearized scheme
\cite{LiaoYanZhang:2018} for a semilinear subdiffusion equation. Thus, the
constant $\Lambda$ is always not very large and the maximum time-step
restriction $\max_{1\le n\le N}\tau_n
\le1/\sqrt[\alpha]{2\pi_A\Gamma(2-\alpha)\Lambda}$
is also not stringent in practical applications.
\end{remark}
\begin{remark}\label{rem: Mittag-Leffler}
The Mittag--Leffler function~$E_\alpha$ also arises naturally in other discrete
and continuous Gr\"{o}nwall inequalities for fractional diffusion and wave
equations \cite[Lemma~2]{Alikhanov2010}, and for weakly singular Volterra
equations~\cite[Theorems 1.3~and 1.6]{DixonMcKee1986}.
The presence of the nonuniform mesh parameter $\rho$ in the argument
of~$E_\alpha$ indicates that sudden, drastic reductions of the time-step should
be avoided.  Nevertheless, %our assumptions still permit local refinement around any desired choice of~$t$.
our discrete Gr\"onwall inequality does not restrict the heterogeneous degree of time mesh, this is, fits for general nonuniform mesh.
\end{remark}

We also have an alternative version of the above theorem.

\begin{theorem}\label{thm: gronwall2}
\cref{thm: gronwall} remains valid if the
condition~\eqref{eq: first Gronwall} is replaced by
\begin{equation}
\sum_{k=1}^nA^{(n)}_{n-k}\diff v^k\le \sum_{k=1}^n\lambda_{n-k}v^{k-\theta}+g^n
    \quad\text{for $1\le n\le N$.}
\end{equation}
Moreover, if the given sequence $(\lambda_l)_{l=0}^{N-1}$ is non-positive and the constant $\Lambda\le0$,
\begin{equation}
v^n\le v^0+\sum_{j=1}^n P^{(n)}_{n-j}g^j\le v^0+\pi_A\Gamma(1-\alpha)\max_{1\le j\le n}\{t_j^{\alpha}g^j\}
	\quad\text{for $1\le n\le N$.}
\end{equation}
\end{theorem}
\begin{proof}
The structure of proof is as before.  However, instead of
\eqref{eq: gronwall 1}~and \eqref{eq: gronwall 2}, we now have
\[
\sum_{j=1}^nP^{(n)}_{n-j}\sum_{k=1}^jA^{(j)}_{j-k}\diff v^k
    \le
    \sum_{j=1}^nP^{(n)}_{n-j}\sum_{k=1}^j\lambda_{j-k}v^{k-\theta}
    +\sum_{j=1}^nP^{(n)}_{n-j}g^j
\]
and
\[
\sum_{j=1}^nP^{(n)}_{n-j}\sum_{k=1}^jA^{(j)}_{j-k}\diff v^k
    =\sum_{k=1}^n\diff v^k\sum_{j=k}^nP^{(n)}_{n-j}A^{(j)}_{j-k}
    =\sum_{k=1}^n\diff v^k=v^n-v^0,
\]
respectively, so that instead of~\eqref{eq: vn2} we obtain
\[
v^n\le v^0
    +\sum_{j=1}^nP^{(n)}_{n-j}\sum_{k=1}^j\lambda_{j-k}v^{k-\theta}
    +\sum_{j=1}^nP^{(n)}_{n-j}g^j.
\]
As before, if $v^1\le v^0$ then $v^1\le G_1$.  For the alternative
case~$v^1>v^0$, we again have $v^{1-\theta}\le v^1$ which now yields
\[
v^1\le v^0+P^{(1)}_0g^1+P^{(1)}_0\lambda_0v^{1-\theta}
    =G_1+P^{(1)}_0\lambda_0v^{1-\theta}\le G_1+\tfrac12 v^1,
\]
where the final step again relies on the step size assumption to
ensure~\eqref{eq: P lambda}.  Thus, once again, $v^1\le2G_1$.  In the inductive
step, \eqref{eq: gronwall 3} is replaced by
\[
v^n\le v^0+\sum_{j=1}^n P^{(n)}_{n-j}g^j
    +\sum_{j=1}^{n-1}P^{(n)}_{n-j}\sum_{k=1}^j\lambda_{j-k}v^{k-\theta}
    +v^nP^{(n)}_0\sum_{k=1}^n\lambda_{n-k},
\]
and by again using~\eqref{eq: gronwall 4} together with the limitation on the
maximum step size, we see that
\[
v^n\le\biggl(
G_n+\sum_{j=1}^{n-1}P^{(n)}_{n-j}\sum_{k=1}^j\lambda_{j-k}v^{k-\theta}\biggr)
        +\frac{v^n}{2},
\]
which is equivalent to~\eqref{eq: gronwall 5} so the remainder of the proof is
unchanged.
\end{proof}

\begin{remark}\label{rem: Jin gronwall}
The discrete fractional Gr\"{o}nwall inequalities in
\cref{thm: gronwall,thm: gronwall2} are valid on very general nonuniform time meshes
and differ substantially from the discrete fractional Gr\"{o}nwall inequality
of Jin et al. \cite[Theorem~2.8]{JinLiZhou2018}, which is built on the uniform mesh for both the L1 scheme
and the convolution quadratures generated by backward difference formulas.
\end{remark}

\begin{remark}[Multi-term and distributed-order Caputo derivatives]\label{rem: other Caputo}
Note that our theory starts only from the discrete convolution form
\eqref{eq: discrete Caputo} and the three assumptions \Ass{1}--\Ass{3},
but not the continuous counterpart \eqref{eq: Caputo}.
Correspondingly, the complementary discrete kernels~$P^{(n)}_{n-j}$
defined in \eqref{eq: P weights} are also independent of \eqref{eq: Caputo}.
In other words, the fractional order~$\alpha$ of Caputo's
derivative~$\fd{\alpha} v$ in \cref{lem: Pn A1 A2,lem: A1-A3},
and the fractional exponent~$\alpha$ in the Mittag--Leffler function~$E_\alpha$
in \cref{lem: E_alpha,thm: gronwall,thm: gronwall2}, are
determined only by the integrand function $\omega_{1-\alpha}(t_n-s)$ of the
lower bound in \Ass{2}, but are independent of
the continuous counterpart of \eqref{eq: discrete Caputo}.

To explain this point more clearly, suppose that the discrete convolution
form \eqref{eq: discrete Caputo} arises from some numerical formula
for a multi-term Caputo derivative
$\sum_{i=1}^{m}w_i\fd{\alpha_i} v$ with $0<\alpha_i<1$ and the weights $w_i>0$,
see \cite{Podlubny:1999}. Then all of the fractional exponents $\alpha_i$ or
the maximum order $\max_{1\leq i\leq m}\alpha_i$ can
determine a single fractional exponent $\alpha$ for \Ass{2} and
the Mittag--Leffler function~$E_\alpha$ in \cref{thm: gronwall,thm: gronwall2}.
Hence, the presented results would be also useful for studying numerical
approximations of multi-term Caputo derivatives and distributed-order Caputo
derivatives, since the latter can be approximated by certain multi-term
derivatives via a proper quadrature rule~\cite{LiaoLyuVongZhao:2017}.
\end{remark}

\begin{remark}[Caputo BDF2-like formula and an open problem]
\label{rem: Caputo BDF2}
There are other practically important formulas, such as the Caputo BDF2-like
approach \cite{GaoSunZhang:2014,LiaoLyuVongZhao:2017,LvXu:2016}.
To start the time-stepping process, one computes the first-level solution
by the L1 approach in \cref{exam: L1}, $(\dfd{\alpha}v)^{1}:=a^{(1)}_{0}\diff
v^1$,
or the Alikhanov formula in \cref{exam: Alikhanov},
$(\dfd{\alpha}v)^{1}:=\hat{a}^{(1)}_{0}\diff v^1$.
For any time-level $t_n$ with $n\ge2$, taking $\theta=0$ and
applying the quadratic polynomial interpolation $\Pi_{2,k}v$,
we have a Caputo BDF2-like formula \cite{LvXu:2016}
\begin{align*}
(\dfd{\alpha}v)^{n}
	\defeq&\,\sum_{k=1}^{n-1}\int_{t_{k-1}}^{t_k}\omega_{1-\alpha}(t_{n}-s)\bra{\Pi_{2,k}v}'(s)\zd{s}\\
&\,+\int_{t_{n-1}}^{t_n}\omega_{1-\alpha}(t_{n}-s)\bra{\Pi_{2,n-1}v}'(s)\zd{s}\quad\text{for $n\ge2$}.
\end{align*}
One can obtain the compact form \eqref{eq: discrete Caputo}
with the discrete kernels $A_{n-k}^{(n)}$,
\begin{equation*}
A^{(n)}_{n-k}\defeq\begin{cases}
	a^{(n)}_0+\rho_{n-1}\brab{b^{(n)}_1+b^{(n)}_0},
	&\text{for $k=n$,}\vspace{0.1cm}\\
	a^{(n)}_{1}+\rho_{n-2}b^{(n)}_2-\brab{b^{(n)}_{1}+b^{(n)}_{0}},
	&\text{for $k=n-1$,}\vspace{0.1cm}\\
    a^{(n)}_{n-k}+\rho_{k-1}b^{(n)}_{n-k+1}-b^{(n)}_{n-k},
	&\text{for $2\le k\le n-2$,}\vspace{0.1cm}\\
	a^{(n)}_{n-1}-b^{(n)}_{n-1},
	&\text{for $k=1$.}
\end{cases}
\end{equation*}
where the coefficients $a_{n-k}^{(n)}$ are defined in \cref{exam: L1},
and the $b_{n-k}^{(n)}$ are defined by
\begin{align*}
&b^{(n)}_{0}\defeq\frac{2}{\tau_{n-1}(\tau_{n-1}+\tau_{n})}\int_{t_{n-1}}^{t_n}
	(s-t_{n-\frac12})\omega_{1-\alpha}(t_{n}-s)\zd{s},\\
&b^{(n)}_{n-k}\defeq\frac{2}{\tau_{k}(\tau_{k}+\tau_{k+1})}\int_{t_{k-1}}^{t_k}
	(s-t_{k-\frac12})\omega_{1-\alpha}(t_{n}-s)\zd{s}\quad \text{for $1\le k\leq n-1$}.
\end{align*}

Notice that if the fractional order $\alpha\to1$, then
$\omega_{3-\alpha}(t)\to t$, $\omega_{2-\alpha}(t)\to1$
and $\omega_{1-\alpha}(t)\to0$, uniformly for~$t>0$.  Thus, we have
$a^{(n)}_0=\omega_{2-\alpha}(\tau_n)/\tau_n\to1/\tau_n$ and
$$b^{(n)}_0=\frac{2}{\tau_{n-1}(\tau_{n-1}+\tau_{n})}
\kbra{\omega_{3-\alpha}(\tau_n)-\frac{\tau_n}{2}\omega_{2-\alpha}(\tau_n)}\to
\frac{\tau_n}{\tau_{n-1}(\tau_{n-1}+\tau_{n})},$$
whereas $a^{(n)}_{n-k}\to0$~and $b^{(n)}_{n-k}\to0$ for $k\ge 1$.
So, when the fractional order $\alpha\to1$,
\begin{align*}
(\dfd{\alpha}v)^{n}\to D_2v^n:=\braB{\frac{1}{\tau_n}+\frac{1}{\tau_{n-1}+\tau_{n}}}\diff v^n
-\frac{\tau_n}{\tau_{n-1}(\tau_{n-1}+\tau_{n})}\diff v^{n-1}
\end{align*}
which is the second-order BDF2 scheme for the classical diffusion equations.
We see that the second kernel $A_1^{(n)}$ can be negative, at least,
when $\alpha$ is close to $1$ (whereas the Caupto BDF2 scheme is shown in \cite{LiaoLyuVongZhao:2017}
to preserve the discrete maximum principle and nonnegativity property when
$\alpha$ is close to 0).

The Caputo BDF2 formula may not meet
our \emph{a priori} assumptions \Ass{1}--\Ass{2}, which results in that our Gr\"{o}nwall inequality
would be not applicable directly. It is not surprising because, for a classical
parabolic equation, the standard discrete Gr\"{o}nwall inequality can also not
be applied directly to the second-order BDF2 scheme.
However, a weighted recombination technique works well;
see the detailed analysis by Thom\'ee~\cite[Theorem 1.7]{Thomee:2006}
for a uniform time mesh, and a similar technique for nonuniform meshes
\cite{Becker:1998,Emmrich:2005}. For the Caputo BDF2 formula, \cref{thm:
gronwall,thm: gronwall2} would be also useful for the stability and convergence
analysis if it can be rearranged to meet the positive and monotone assumptions
\Ass{1}--\Ass{2}. On the uniform mesh with $\tau_n=\tau$, Lv and Xu
\cite{LvXu:2016} developed a new technique of variable-weights recombination and
achieved a new form of $(\dfd{\alpha}v)^{n}$ with a new variable
$\bar{v}^k:=v^k-\eta v^{k-1}$ and~$\bar{v}^0:=v^0$; in our notations,
\begin{equation}\label{eq: Caputo BDF2 new}
(\dfd{\alpha}v)^{n}=\sum_{k=1}^n
	\bar{A}^{(n)}_{n-k}\diff \bar{v}^k+v^0\sum_{j=1}^nA^{(n)}_{n-j}\eta^j
\end{equation}
where the combination parameter $\eta:=\frac12\brab{1-A^{(n)}_{1}/A^{(n)}_{0}}$.
From the substitution formulas
\begin{align*}
v^k%%=\bar{v}^k+\eta v^{k-1}=\bar{v}^k+\eta\brat{\bar{v}^{k-1}+\eta v^{k-2}}=\cdots
=\sum_{\ell=0}^{k}\eta^{k-\ell}\bar{v}^{\ell}\quad\text{and}\quad
\diff v^k=\sum_{\ell=1}^{k}\eta^{k-\ell}\diff\bar{v}^{\ell}+\eta^kv^0,
\end{align*}
one has a new series of discrete convolution weights
\begin{equation*}
\bar{A}^{(n)}_{n-k}:=\sum_{j=k}^{n}A^{(n)}_{n-j}\eta^{j-k}\quad\text{for $1\leq k\leq n$}.
\end{equation*}
The results of \cite[Lemma 3.2]{LvXu:2016} imply that $0<\eta<2/3$
and the new convolution kernels $\bar{A}^{(n)}_{n-k}$ are positive and monotone,
\begin{equation*}
\bar{A}^{(n)}_{0}>\bar{A}^{(n)}_{1}>\cdots>\bar{A}^{(n)}_{n-1}>0\quad\text{for $1\leq k\leq n$}.
\end{equation*}
Thus, our discrete Gr\"{o}nwall inequalities
 (and the complementary discrete convolution kernels~$P^{(n)}_{n-j}$ as well)
could be applied  for this new form \eqref{eq: Caputo BDF2 new}
directly once a proper constant $\pi_A$ in \Ass{2} is determined by a more careful examination.

Nonetheless, we do not know whether the variable-weights recombination
technique~\cite{LvXu:2016} works on nonuniform time grids. More precisely, it
has yet to be determined what constraints must be imposed on a nonuniform
mesh so that the new discrete form \eqref{eq: Caputo BDF2 new} satisfies
the \emph{a priori} assumptions \Ass{1}-\Ass{3} required by
\cref{thm: gronwall,thm: gronwall2}. This problem could be very challenging, at
least technically, and remains open to us.
\end{remark}

%%%%%%%%%%%%%%%%%%%%%%%%%%%%%%%%%%%%%%%%%%%%%%%%%%%%%%%%%%%%%%%%%%%%%%%%%%%%%%%
\section{Stability and consistency}\label{sec: stability consistency}
We will now outline how the results of Section~\ref{sec: discrete Gronwall}
can be applied to study a numerical solution of problem~\eqref{eq: reaction-diffusion}.
For simplicity, we restrict our attention to the case of a linear reaction term
$f(x,t,u)\defeq\kappa u+\psi(x,t)$ with a constant~$\kappa\ge0$.
By applying the first Green identity, the
fractional PDE~\eqref{eq: reaction-diffusion} is written in a weak form as
\begin{equation}\label{eq: weak}
\iprod{\fd{\alpha} u,v}+\Bilin(u,v)=\kappa\iprod{u,v}+\iprod{\psi(t),v}
	\quad\text{for all $v\in H^1_0(\Omega)$ and for $0<t\le T$,}
\end{equation}
where $\iprod{u,v}$ denotes the inner product in~$L_2(\Omega)$, and
$\Bilin(u,v)=\iprod{\opL u,v}$ is the bilinear form induced by the elliptic
operator~$\opL$. Since the latter is strongly elliptic, by increasing
$\kappa$ if necessary, we may assume that the bilinear form is coercive: there is a constant $c>0$ such that
\begin{equation}\label{eq: coercive}
\Bilin(v,v)\ge c\|v\|_{H^1_0(\Omega)}^2\quad\text{for all $v\in H^1_0(\Omega)$.}
\end{equation}

Let $X_h$ be a finite dimensional subspace of~$H^1_0(\Omega)$; for example,
a (conforming) finite element space based on a triangulation of~$\Omega$
with the mesh size~$h$.  Galerkin's method yields a spatially-discrete approximate
solution~$u_h:[0,T]\to X_h$ satisfying
\begin{equation}\label{eq: spatially discrete}
\iprodb{\fd{\alpha} u_h,\chi}+\Bilin(u_h,\chi)=\kappa\iprodb{u_h,\chi}
	+\iprodb{\psi(t),\chi}
	\quad\text{for all $\chi\in X_h$ and $0<t\le T$,}
\end{equation}
with $u_h(0)=u_{h0}\approx u_0$ for a suitable $u_{h0}\in X_h$.
To compute a fully-discrete numerical solution~$u^n_h\in X_h$, where
$u(t_n)\approx u^n_h$ for~$1\le n\le N$, we apply the approximation~\eqref{eq:
discrete Caputo} to the fractional derivative term
in~\eqref{eq: spatially discrete} so that
\begin{equation}\label{eq: fully discrete}
\iprodb{(\dfd{\alpha}u_h)^{n-\theta},\chi}+\Bilin\brab{u_h^{n-\theta},\chi}
	=\kappa\iprodb{u_h^{n-\theta},\chi}+\iprodb{\psi(t_{n-\theta}), \chi }
\end{equation}
for all $\chi\in X_h$ and for~$1\le n\le N$.

The next lemma is a
discrete analogue of the inequality~\cite[Lemma~1]{Alikhanov2010}
\[
\brab{\fd{\alpha} \|v\|^2}(t)\le2\iprodb{(\fd{\alpha} v)(t),v(t)}
    \quad\text{for $0\le t\le T$ and $0<\alpha<1$,}
\]
and helps set the stage for applying our discrete fractional Gr\"{o}nwall
inequality.

\begin{lemma}\label{lem: energy 2}
Let the assumption \Ass{1} hold and  fix the parameter~$\theta\in[0,1)$.  Then
every sequence~$(v^n)_{n=0}^N$ in~$L_2(\Omega)$ satisfies
\[
\sum_{k=1}^nA^{(n)}_{n-k}\diff\brab{\mynormt{v^k}^2}
	\le 2\iprodb{(\dfd{\alpha}v)^{n-\theta},v^{n-\theta}}-
	d_n\brab{\theta^{(n)}-\theta}\mynormb{(\dfd{\alpha}v)^{n-\theta}}^2,
\]
for $1\le n\le N$, where $0<d_n<1/A^{(n)}_0$ and $0<\theta^{(n)}<1/2$ are given
by
\[
d_n\defeq\frac{2A^{(n)}_0-A^{(n)}_1}{A^{(n)}_0(A^{(n)}_0-A^{(n)}_1)}>0
\quad\text{and}\quad
\theta^{(n)}\defeq\frac{A^{(n)}_0-A^{(n)}_1}{2A^{(n)}_0-A^{(n)}_1}<\frac{1}{2}.
\]
\end{lemma}
\begin{proof}
By \cref{lem: energy 1} (see appendix~\ref{sec: technical}),
\begin{multline*}
2\iprodb{(\dfd{\alpha}v)^{n-\theta},v^{n-\theta}}
    =2\theta\iprodb{\dfd{\alpha}v)^{n-\theta},v^{n-1}}
    +2(1-\theta)\iprodb{\dfd{\alpha}v)^{n-\theta},v^n}\\
	\ge\sum_{k=1}^nA^{(n)}_{n-k}\bra{\mynormb{v^k}^2-\mynormb{v^{k-1}}^2}+
	\biggl(\frac{1-\theta}{A^{(n)}_0}-\frac{\theta}{A^{(n)}_0-A^{(n)}_1}\biggr)
	\mynormb{(\dfd{\alpha}v)^{n-\theta}}^2,
\end{multline*}
and the second term on the right side equals
$d_n(\theta^{(n)}-\theta)\|(\dfd{\alpha}v)^{n-\theta}\|^2$.
\end{proof}

\begin{theorem}\label{thm: stability}
Let the assumption~\Ass{1} hold and $0\le\theta\le\theta^{(n)}$
for~$1\le n\le N$.  Then the fully-discrete solution~$u^n_h\in X_h$, defined
by~\eqref{eq: fully discrete}, satisfies
\[
\sum_{k=1}^nA^{(n)}_{n-k}\diff\brab{\mynormb{u^k_h}^2}
	\le2\kappa\mynormb{u^{n-\theta}_h}^2
	+2\mynormb{u^{n-\theta}_h}\mynormb{\psi(t_{n-\theta})}
	\quad\text{for $1\le n\le N$.}
\]
\end{theorem}
\begin{proof}
Put $\chi=2u_h^{n-\theta}$ in the Galerkin discrete equation~\eqref{eq: fully discrete},
apply \cref{lem: energy 2} with~$v^n=u^n_h$, and use
positive-definiteness~\eqref{eq: coercive} of the bilinear form.
\end{proof}

Applying the discrete fractional Gr\"{o}nwall inequality from
\cref{thm: gronwall} with
\[
v^n\defeq\mynormb{u^n_h},\quad
g^n\defeq2\mynormb{\psi(t_{n-\theta})},\quad\lambda_0\defeq 2\kappa\quad
\text{and}\quad\text{$\lambda_j\defeq0$ for $1\le j\le N-1$,}
\]
we see from \cref{thm: stability} that the
scheme~\eqref{eq: fully discrete} is stable in $L_2(\Omega)$,
\[
\mynormb{u^n_h}\le2E_\alpha\bigl(4\max(1,\rhomax)\pi_A\kappa t_n^\alpha\bigr)
\biggl(\mynormb{u_{0h}}
	+2\max_{1\le k\le n}
	\sum_{j=1}^kP^{(k)}_{k-j}\mynormb{\psi(t_{j-\theta})}\biggr),
\]
provided $\theta\le\theta^{(n)}$ for~$1\le n\le N$. The inequality from
Remark~\ref{rem: sum P} yields a weaker but simpler stability estimate,
\[
\mynormb{u^n_h}\le2E_\alpha\bigl(4\max(1,\rhomax)\pi_A\kappa t_n^\alpha\bigr)
\biggl(\mynormb{u_{0h}}+2\pi_A\Gamma(1-\alpha)\max_{1\le k\le n}t_k^{\alpha}\mynormb{\psi(t_{k-\theta})}\biggr).
\]

To bound the error in~$u^n_h$, we introduce the Ritz
projector~$R_h:H^1_0(\Omega)\to X_h$, which is well-defined by
\begin{equation}\label{eq: Ritz}
\Bilin(R_hv,\chi)=\Bilin(v,\chi)
	\quad\text{for all $v\in H^1_0(\Omega)$ and $\chi\in X_h$,}
\end{equation}
because the bilinear form satisfies~\eqref{eq: coercive}.
Put $e^n_h=u^n_h-R_hu^n\in X_h$ where $u^n=u(t_n)$, so that
\[
\mynormb{u^n_h-u^n}\le\mynormb{u^n-R_hu^n}+\mynormb{e^n_h}.
\]
The error in the Ritz projection~$R_hu^n$ is estimated in the usual way from a
study of the elliptic problem, so it suffices to deal with~$\mynormb{e^n_h}$.
Using the weak form \eqref{eq: weak} at $t=t_{n-\theta}$,
with~$v=\chi$, we see that
\begin{equation}\label{eq: weak n-theta}
\iprodb{(\fd{\alpha} u)(t_{n-\theta}),\chi}+\Bilin(u(t_{n-\theta}),\chi)
    =\kappa\iprodb{u(t_{n-\theta}),\chi}+\iprodb{\psi(t_{n-\theta}),\chi}.
\end{equation}
%where $u^{n-\theta}\defeq\theta u^{n-1}+(1-\theta)u^n$.
%\[
%\brab{\fd{\alpha} u}^{n-\theta}=\theta\brab{\fd{\alpha} u}(t_{n-1})
%    +(1-\theta)\brab{\fd{\alpha} u}(t_n)
%\quad\text{and}\quad
%u^{n-\theta}\defeq\theta u^{n-1}+(1-\theta)u^n.
%\]
It follows from~\eqref{eq: fully discrete} that
\begin{align*}
\iprodb{(\dfd{\alpha}e_h)^{n-\theta},\chi}+\Bilin(e_h^{n-\theta},\chi)
    &=\kappa\iprodb{u_h^{n-\theta},\chi}+\iprodb{\psi(t_{n-\theta}),\chi}\\
    &\qquad{}-\iprodb{(\dfd{\alpha}R_hu)^{n-\theta},\chi}
    -\Bilin(R_hu^{n-\theta},\chi).
\end{align*}
Therefore, since \eqref{eq: Ritz}~and \eqref{eq: weak n-theta} imply
\begin{align*}
\Bilin(R_hu^{n-\theta},\chi)=&\,\Bilin(u^{n-\theta}-u(t_{n-\theta}),\chi)
	+\Bilin(u(t_{n-\theta}),\chi)\\
    =&\,-\iprod{\triangle(u^{n-\theta}-u(t_{n-\theta})),\chi}
	+\kappa\iprodb{u(t_{n-\theta}),\chi}\\&\,+\iprodb{\psi(t_{n-\theta}),\chi}
        -\iprodb{(\fd{\alpha} u)(t_{n-\theta}),\chi},
\end{align*}
we have
\[
\iprodb{(\dfd{\alpha}e_h)^{n-\theta},\chi}+\Bilin(e_h^{n-\theta},\chi)
	=\kappa\iprodb{e_h^{n-\theta},\chi}+\iprod{\mathcal{R}^n,\chi}
\quad\text{for all $\chi\in X_h$,}
\]
where
\[
\mathcal{R}^n=
	(\fd{\alpha} u)(t_{n-\theta})-(\dfd{\alpha}R_hu)^{n-\theta}
	-\kappa\bigl(u(t_{n-\theta})-R_hu^{n-\theta}\bigr)
	+\triangle\bigl(u^{n-\theta}-u(t_{n-\theta})\bigr).
\]
Choosing $\chi=2e^{n-\theta}_h$ and arguing as before,
but now with $v^n:=\mynormb{e^n_h}$~and
$g^n:=2\mynormb{\mathcal{R}^n}$, we see that (for appropriate $\theta$)
\[
\mynormb{e^n_h}\le2E_\alpha\bigl(4\max(1,\rhomax)\pi_A\kappa t_n^\alpha\bigr)
	\biggl(\|u_{0h}-u_0\|
    +2\max_{1\le k\le n}\sum_{j=1}^kP^{(k)}_{k-j}\|\mathcal{R}^j\|\biggr)
\]
for $1\le n\le N$.  A complete error analysis would typically proceed by
applying the triangle inequality to obtain
\begin{align*}
\mynormb{\mathcal{R}^j}
	&\le\mynormb{(\fd{\alpha} u)(t_{j-\theta})-(\dfd{\alpha}u)^{j-\theta}}
	+\mynormb{(\dfd{\alpha}(u-R_hu))^{j-\theta}}\\
	&\qquad{}+\kappa\bigl\|(u-R_hu)^{j-\theta}\bigr\|
	+\bigl\|(\kappa+\triangle)\bigl(u^{j-\theta}-u(t_{j-\theta})\bigr)\bigr\|,
\end{align*}
and estimating separately the resulting convolutional sums over~$j$,
refer to a new technique of global consistency error analysis
developed in recent works~\cite{LiaoLiZhang:2018,LiaoMcLeanZhang:2018,LiaoYanZhang:2018}.
The details would depend on the choice of the discrete kernels $A^{(n)}_{n-j}$~and
of the space~$X_h$, and would rely on some \emph{a priori} estimates for the partial derivatives of~$u$.

A similar approach
works if finite differences are used for the space
discretization~\cite{LiaoLiZhang:2018}, by introducing an appropriate discrete
$\ell_2$ inner product in place of the inner product $\iprod{u,v}$.
%%%%%%%%%%%%%%%%%%%%%%%%%%%%%%%%%%%%%%%%%%%%%%%%%%%%%%%%%%%%%%%%%%%%%%%%%%%%%%%

\section*{Acknowledgements}
Hong-lin Liao and Jiwei Zhang would like to thank Prof. Ying Zhao,  Prof. Weiwei Sun,
Prof. Martin Stynes and Dr. Yonggui Yan for their valuable discussions and fruitful suggestions.
Hong-lin Liao thanks for the hospitality of Beijing
Computational Science Research Center during the period of his visit.
%We also thank the anonymous referees for their valuable comments
%and suggestions, which are very helpful for improving the quality of the paper.

\appendix
\section{Two technical inequalities}\label{sec: technical}
The proof of \cref{lem: energy 2} relies on the following result,
essentially due to Alikhanov~\cite[Lemma~1]{Alikhanov2015}.

\begin{lemma}\label{lem: energy 1}
If the assumption \Ass{1} holds, then every sequence~$(v^n)_{n=0}^N$
in~$L_2(\Omega)$ satisfies
\[
2\iprod{(\dfd{\alpha}v)^{n-\theta},v^n}
\smash{\ge\sum_{k=1}^nA^{(n)}_{n-k}\braB{\mynormb{v^k}^2-\mynormb{v^{k-1}}^2}
	+\frac{\mynormb{(\dfd{\alpha}v)^{n-\theta}}^2}{A^{(n)}_0}}
\]
and
\[
2\iprod{(\dfd{\alpha}v)^{n-\theta},v^{n-1}}
\smash[t]{\ge\sum_{k=1}^nA^{(n)}_{n-k}\braB{\mynormb{v^k}^2-\mynormb{v^{k-1}}^2}
	-\frac{\mynormb{(\dfd{\alpha}v)^{n-\theta}}^2}{A^{(n)}_0-A^{(n)}_1}}
\]
for $1\le n\le N$, provided we set $A^{(1)}_1=0$ in the case~$n=1$.
\end{lemma}
\begin{proof}
Fix~$n$ and consider the difference
\[
J_n:=2\iprod{(\dfd{\alpha}v)^{n-\theta},v^n}
	-\sum_{k=1}^nA^{(n)}_{n-k}\braB{\mynormb{v^k}^2-\mynormb{v^{k-1}}^2}.
\]
We have
\begin{align*}
J_n&=\sum_{k=1}^nA^{(n)}_{n-k}\left(2\iprod{v^k-v^{k-1},v^n}
		-\iprod{v^k-v^{k-1},v^k+v^{k-1}}\right)\\
	&=\sum_{k=1}^nA^{(n)}_{n-k}\iprod{v^k-v^{k-1},2v^n-(v^k+v^{k-1})}
\end{align*}
and, using the identity
$2v^n-(v^k+v^{k-1})=v^k-v^{k-1}+2\sum_{j=k+1}^n(v^j-v^{j-1})$,
\begin{align*}
J_n&=\sum_{k=1}^nA^{(n)}_{n-k}\mynormb{v^k-v^{k-1}}^2
	+2\sum_{k=1}^nA^{(n)}_{n-k}\sum_{j=k+1}^n
		\iprod{v^k-v^{k-1},v^j-v^{j-1}}\\
	&=\sum_{k=1}^nA^{(n)}_{n-k}\mynormb{v^k-v^{k-1}}^2
	+2\sum_{j=2}^n\sum_{k=1}^{j-1}A^{(n)}_{n-k}
		\iprod{v^k-v^{k-1},v^j-v^{j-1}}.
\end{align*}
To continue the proof, it is convenient to introduce
\[
w^j:=\sum_{k=1}^jA^{(n)}_{n-k}(v^k-v^{k-1})
	\quad\text{and}\quad Q_j:=\frac{1}{A^{(n)}_{n-j}}\quad
	\text{for~$1\le j\le n$.}
\]
Notice that $v^j-v^{j-1}=Q_j(w^j-w^{j-1})$ for~$2\le j\le n$, and that
the assumption \Ass{1} implies $Q_1\ge Q_2\ge\cdots\ge Q_n$.
Thus, one deduces that
\begin{align*}
J_n&=Q_1\mynormb{w^1}^2+\sum_{j=2}^nQ_j\mynormb{w^j-w^{j-1}}^2
	+2\sum_{j=2}^nQ_j\iprod{w^{j-1},w^j-w^{j-1}}\\
	&=Q_1\mynormb{w^1}^2+\sum_{j=2}^nQ_j\braB{\mynormb{w^j}^2
        -\mynormb{w^{j-1}}^2}\\
	&=Q_n\mynormb{w^n}^2+\sum_{j=1}^{n-1}(Q_j-Q_{j+1})
        \mynormb{w^j}^2\ge Q_n\mynormb{w^n}^2.
\end{align*}
The first inequality now follows by noting that $w^n=(\dfd{\alpha}v)^{n-\theta}$~and
$Q_n=1/A^{(n)}_0$. Furthermore, by using the identity
$v^{n-1}=v^n-(v^n-v^{n-1})=v^n-Q_n(w^n-w^{n-1})$, we have
\begin{align*}
2\langle(\dfd{\alpha}v)^{n-\theta},&v^{n-1}\rangle
	-\sum_{k=1}^nA^{(n)}_{n-k}\braB{\mynormb{v^k}^2-\mynormb{v^{k-1}}^2}
	=J_n-2Q_n\iprod{w^n,w^n-w^{n-1}}\\
	&\ge Q_n\mynormb{w^n}^2
        +(Q_{n-1}-Q_n)\mynormb{w^{n-1}}^2-2Q_n\iprod{w^n,w^n-w^{n-1}}\\
	&=-Q_n\mynormb{w^n}^2+2Q_n\iprod{w^n,w^{n-1}}
        +(Q_{n-1}-Q_n)\mynormb{w^{n-1}}^2\\
	&=\frac{1}{Q_{n-1}-Q_n}\braB{
		\mynorm{Q_nw^n+(Q_{n-1}-Q_n)\,w^{n-1}}^2
		-Q_nQ_{n-1}\mynorm{w^n}^2}\\
    &\geq-\frac{Q_nQ_{n-1}}{Q_{n-1}-Q_n}\mynorm{w^n}^2
    =\frac{-\mynorm{w^n}^2}{A^{(n)}_0-A^{(n)}_1}.
\end{align*}
Therefore the claimed second inequality follows and the proof is complete.
\end{proof}
%%%%%%%%%%%%%%%%%%%%%%%%%%%%%%%%%%%%%%%%%%%%%%%%%%%%%%%%%%%%%%%%%%%%%%%%%%%%%%%
%\bibliographystyle{siamplain}
%\bibliography{LiaoEtAl_Gronwallrefs}

%%%%%%%%%%%%%%%%%%%%%%%%%%%%%%%%%%%%%%%%%%%%%%%%%%%%%%%%%%%%%%%%%%%%%%%%%%%%%%%
\end{document}